\documentclass[12pt]{article}
\usepackage{amsfonts,url}
\usepackage{authblk}

\newtheorem{theorem}{Theorem}[section]
\newtheorem{prop}[theorem]{Proposition}

\newtheorem{conj}{Conjecture}[section]

\newtheorem{unnumber}{}

\newtheorem{prepf}[unnumber]{Proof}
\newenvironment{pf}{\prepf\rm}{\endprepf}

\begin{document}

\title{Algebraic properties of chromatic roots}

\author[1]{Peter J. Cameron\thanks{pjc20@st-andrews.ac.uk}}
\author[2]{Kerri Morgan\thanks{Kerri.Morgan@monash.edu}}
\affil[1]{School of Mathematics and Statistics, University of St Andrews}
\affil[2]{Faculty of Information Technology, Monash University}
\date{}
\maketitle

\begin{abstract}
A \emph{chromatic root} is a root of the chromatic polynomial of a graph.
Any chromatic root is an algebraic integer. Much is known about
the location of chromatic roots in the real and complex numbers, but
rather less about their properties as algebraic numbers. This question
was the subject of a seminar at the Isaac Newton Institute in late 2008.
The purpose of this paper is to report on the seminar and subsequent
developments.

We conjecture that, for every algebraic integer $\alpha$, there is a
natural number $n$ such that $\alpha+n$ is a chromatic root. This is proved
for quadratic integers; an extension to cubic integers has been found by
Adam Bohn. The idea is to consider certain special classes of graphs
for which the chromatic polynomial is a product
of linear factors and one ``interesting'' factor of larger degree. We also
report computational results on the Galois groups of irreducible factors
of the chromatic polynomial for some special graphs. Finally,
extensions to the Tutte polynomial are mentioned briefly.
\end{abstract}

\section{Chromatic roots}
\label{s1}

A proper \emph{colouring} of a graph $G$ is a function from the vertices
of $G$ to a set of $q$ colours with the property that adjacent vertices
receive different colours. The \emph{chromatic polynomial} $P_G(q)$ of $G$
is the function whose value at the positive integer $q$ is the number of
proper colourings of $G$ with $q$ colours. It is a monic polynomial in $q$
with integer coefficients, whose degree is the number of vertices of $G$.

A \emph{chromatic root} is a complex number $\alpha$ which is a root of some
chromatic polynomial.

\subsection{Location of chromatic roots}

A lot of attention has been paid to the location of chromatic roots in the
complex plane, which we now outline.

\paragraph{Integer quadratic roots}

An integer $m$ is a root of $P_G(q) = 0$ if and only if the chromatic number
of $G$ (the smallest number of colours required for a proper colouring of $G$)
is greater than $m$. Hence every non-negative integer is a chromatic root. (For
example, the complete graph $K_{m+1}$ cannot be coloured with $m$ colours.)

On the other hand, no negative integer is a chromatic root.

\paragraph{Real chromatic roots}

The non-trivial parts of the following theorem are due to Jackson~\cite{r9} and
Thomassen~\cite{r18}.

\begin{theorem}
\begin{enumerate}
\item There are no negative chromatic roots, none in the interval $(0, 1)$,
and none in the interval $(1, \frac{32}{27}]$.
\item
Chromatic roots are dense in the interval $[\frac{32}{27},\infty)$.
\end{enumerate}
\end{theorem}

\paragraph{Complex chromatic roots}

For some time it was thought that chromatic roots must have non-negative
real part. This is true for graphs with fewer than ten vertices. But
Sokal~\cite{r16} showed the following.

\begin{theorem}
Complex chromatic roots are dense in the complex plane.
\end{theorem}

This is connected with the Yang--Lee theory of phase transitions. Sokal used
\emph{theta-graphs} in his proof; these graphs also play a role in our
investigation.

\subsection{Algebraic integers}

An \emph{algebraic number} is a (complex) root of a polynomial over the
integers; an \emph{algebraic integer} is a root of a monic integer polynomial.
Clearly any chromatic root is an algebraic integer. Our main question is the
converse:
\begin{quote}
Which algebraic integers are chromatic roots?
\end{quote}

Let $G + K_n$ denote the graph obtained by adding $n$ new vertices to $G$,
joined to one another and to all existing vertices. Then
\[P_{G+K_n}(q) = q(q-1)\cdots(q-n+1)P_G(q-n).\]
We conclude:

\begin{prop}
If $\alpha$ is a chromatic root, then so is $\alpha + n$, for any natural
number $n$.
\label{p:translate}
\end{prop}

However, the set of chromatic roots is far from being a semiring; it is not
closed under either addition or multiplication. This can be seen as follows.
By Sokal's theorem, we can find a chromatic root $\alpha$ with negative real part
and modulus less than $1$. Then neither $\alpha+\overline{\alpha}$ nor $\alpha\cdot\overline{\alpha}$
is a chromatic root: the first is a negative real number, the second lies in
$(0, 1)$.

As partial replacement, here are two conjectures:

\begin{conj}[The $\alpha + n$ conjecture]
Let $\alpha$ be an algebraic integer. Then
there exists a natural number $n$ such that $\alpha + n$ is a chromatic root.
\end{conj}

\begin{conj}[The $n\alpha$ conjecture]
Let $\alpha$ be a chromatic root. Then $n\alpha$
is a chromatic root for any natural number $n$.
\end{conj}

If the $\alpha+n$ conjecture is true, we can ask, for given $\alpha$, what
is the smallest $n$ for which $\alpha + n$ is a chromatic root?

The $\alpha + n$ conjecture can be reformulated as follows. A monic integer
polynomial $f(q)$ of degree $n$ can be transformed by a substitution $x = q+a$
into a unique monic integer polynomial in which the coefficient of $x^{n-1}$
lies between $0$ and $n-1$ (inclusive). We call such a polynomial
\emph{standard}. The conjecture asserts that every standard irreducible monic
integer polynomial is the standard form of a factor of a chromatic polynomial.
This formulation lends itself more readily to computation.

\subsection{An example}

The \emph{golden ratio} $\alpha = (\sqrt{5}-1)/2$ is an algebraic integer,
since it satisfies $\alpha^2 + \alpha - 1 = 0$. It is not a chromatic root,
as it lies in $(0, 1)$.

Also, $\alpha+1$ and $\alpha+2$ are not chromatic roots, since their algebraic
conjugates are negative or in $(0,1)$. There are, however, graphs (for 
example, the truncated icosahedron) which have chromatic roots very close to
$\alpha + 2$, the so-called ``golden root'' \cite{r2}.

We do not know whether $\alpha + 3$ is a chromatic root or not. However, we
will see that $\alpha + 4$ is a chromatic root (the smallest graph having
chromatic root $\alpha + 4$ has eight vertices), and hence so is $\alpha + n$
for any natural number $n\ge4$.

\paragraph{Remark} 
We showed that $\alpha + 1$ and $\alpha + 2$ are not chromatic roots by
showing that they have conjugates in forbidden regions of the real line. Is
there another technique for proving such negative results? Perhaps resolving
the question whether $\alpha + 3$ is a chromatic root would help with this.

\section{Two reductions}

Let $G$ be a graph which is the union of two graphs $G_1$ and $G_2$, whose
intersection is a complete graph of size $k$. Such a graph is called a
\emph{clique-sum} of $G_1$ and $G_2$.

We have
\[P_G(q)=\frac{P_{G_1}(q)P_{G_2}(q)}{q(q-1)\cdots(q-k+1)}.\]
For, having chosen any colouring of $G_1$ with $q$ colours, a proportion
$1/q(q-1)\cdots(q-k+1)$ of the colourings of $G_2$ agree on this intersection.

So if a graph is a clique-sum of smaller graphs, the irreducible factors of
its chromatic polynomial all occur in these smaller graphs. So we need only
consider graphs which cannot be expressed as clique-sums. In particular, we
can take connected graphs.

The argument in Proposition~\ref{p:translate} shows that, furthermore, we may
assume that there is no vertex joined to all others.

\section{Rings of cliques}

Our strategy is to choose certain special graphs for which the chromatic
polynomial can be computed explicitly. The most productive class we found
are the \emph{rings of cliques}, defined as follows:

Let $a_0,\ldots,a_{k-1}$ be positive integers. The graph 
$R(a_0,\ldots,a_{k-1})$ is the disjoint union of complete subgraphs
$C_0,\ldots,C_{k-1}$ with $a_0,\ldots,a_{k-1}$ vertices
respectively, together with all edges from $C_i$ to $C_{i+1}$ for
$i = 0,\ldots, k-1$ (where indices are taken modulo $k$).

The following theorem was proved by Read~\cite{r13}.

\begin{theorem}
The chromatic polynomial of $R(1,a_1,\ldots,a_{k-1})$ is a product of linear
factors and the polynomial
\[\frac{1}{q}\left(\prod_{i=1}^{k-1}(q-a_i)-\prod_{i=1}^{k-1}(-a_i)\right)\]
of degree $k-2$.
\end{theorem}

We call the displayed polynomial the \emph{interesting factor}.

Read later found a more complicated formula for the chromatic polynomial of
arbitrary rings of cliques~\cite{r14}.

In connection with the $n\alpha$ conjecture, we make the following observation:

\begin{prop}
If $\alpha$ is a root of the interesting factor of $R(1,a_1,\ldots,a_{k-1})$,
then $n\alpha$ is a root of the interesting factor of
$R(1,na_1,\ldots,na_{k-1})$.
\end{prop}

\paragraph{Example 1}
The graph $R(1, 1,\ldots, 1)$ (with $k$ entries $1$) is simply an $k$-cycle.
The interesting factor is
\[\frac{(q-1)^{k-1}-(-1)^{k-1}}{q}=\frac{x^{k-1}-(-1)^{k-1}}{x+1},\]
where we have put $x=q-1$. So the roots have the form $1+\omega$, where
\begin{itemize}\itemsep0pt
\item if $k-1$ is odd, then $\omega$ is a $2(k-1)$th root of unity which is
not a $(k-1)$th root and is not $-1$;
\item if $k-1$ is even, then $\omega$ is a $(k-1)$th root of unity which is 
not $-1$.
\end{itemize}

We conclude that the $\alpha+n$ conjecture is true for roots of unity (and
indeed, if $\omega$ is a root of unity, then $\omega+1$ is a chromatic root).

\paragraph{Example 2} The interesting factor for $R(1,1,1,5)$ is $q^2-7q+11$,
which has a root $\alpha + 4$, where $\alpha$ is the golden ratio. This is
the example promised earlier, and is the smallest graph which has a chromatic
root $\alpha + 4$.

\subsection{Quadratic integers}

In this section, we prove the $\alpha+n$ conjecture for quadratic integers.

\begin{theorem}
Let $\alpha$ be an integer in a quadratic number field. Then there is a
natural number $n$ such that $\alpha+n$ is a chromatic root.
\end{theorem}

\begin{pf}
If $\alpha$ is irrational, then the set $\{\alpha+n:n\in\mathbb{Z}\}$ is the
set of all quadratic integers with given discriminant. So it is enough to show
that, for any non-square $d$ congruent to $0$ or $1$ mod~$4$, there is a 
quadratic integer with discriminant $d$ which is a chromatic root.

The interesting factor of $R(1,1,a,b)$ is $x^2-(a+b+1)x+(ab+a+b)$. The
discriminant of this quadratic is
\[(a+b+1)^2-4(ab+a+b)=(a-b)^2-2(a+b)+1.\]
Now $a+b$ and $a-b$ are integers with the same parity. If they are both even,
say $2l$ and $2m$ with $m<l$, we want $d=4m^2-4l+1$. Any number $d$ congruent
to $1$ mod~$4$ is of this form: choose $m$ such that $4m(m-1)>d-1$, and then
$l=m^2-(d-1)/4$. The argument for $d$ congruent to $0$ mod~$4$ is similar.
\end{pf}

Rings of $k$ cliques, one of size $1$, give ``interesting factors'' of degree
$k-2$, whereas only $k-3$ independent parameters are theoretically required to
prove the $\alpha+n$ conjecture for algebraic integers of degree $k-2$. So
it is possible that these graphs would suffice for the purpose. However,
computational evidence in Section~\ref{s6} suggests that this may be difficult.
We have been unable to find such a polynomial with Galois group $C_5$, for
example.

Here is a table of the smallest graphs we found with real quadratic roots of
given discriminant; we give the number of vertices in the graph, the quadratic
factor, and the graph (given by its number in McKay's list of connected
graphs~\cite{r10}, if it is not a ring of cliques).

\begin{center}
\begin{tabular}{|c|c|c|c|}
\hline
Discriminant & Polynomial & No. of vertices & Graph number \\
\hline
$5$ & $x^2-7x+11$ & $8$ & $R(1,1,1,5)$ \\
$8$ & $x^2-6x+7$ & $9$ & $198748$ \\
$12$ & $x^2-8x+13$ & $9$ & $R(1,1,1,6)$ \\
$13$ & $x^2-7x+9$ & $10$ & $10756635$ \\
\hline
\end{tabular}
\end{center}

\section{Bicliques}
\label{s4}

A \emph{biclique} is a graph whose vertex set is the union of two cliques
$C$ and $D$, of sizes $n$ and $m$, say $D=\{w_1,\ldots,w_m\}$. For 
$i=1,\ldots,m$, let $F_i$ be the set of neighbours of $w_i$ in $C$. We denote
this graph by $A(\mathcal{F})$, where $\mathcal{F}=(F_1,\ldots,F_m)$.
(Think of $m$ as fixed and $n$ arbitrary.)

We may assume without loss of generality that:
\begin{itemize}
\item
The union $U$ of all the sets $F_i$ is the whole of $C$. For, if not, then the
graph is a clique-sum: the subgraphs $D\cup U$ and $C$ intersect in the
clique $U$.
\item
The intersection of the sets $F_i$ is empty. For a vertex in this intersection
is joined to every other vertex in the graph.
\end{itemize}

\begin{prop}
The chromatic polynomial of a biclique can be computed in terms of $n$ and the
sizes of the $m$ sets $F_i$ and their intersections.
\end{prop}

\begin{pf}
First, ignore the edges within $D$. If $q$ colours are available, then
$C$ can be coloured in $(q)_n = q(q-1)\cdots(q-n+1)$ ways, and so it is
enough to count the number of colourings of the vertices in $D$; there are
$q-|F_i|$ ways to colour $w_i$, and so the number of colourings is the product
of these numbers. We have to count the subset of these colourings in which
all the vertices in $D$ receive different colours. This can be done by
M\"obius inversion over the poset of partitions of $\{1,\ldots,n\}$ (whose
M\"obius function is known, see \cite{r15}), if we can compute, for each such
partition, the number of colourings in which vertices with indices in the
same part have the same colour. If $I$ is a part, then there are
$q-|\bigcup_{i\in I}F_i|$ ways to choose this colour, and multiplying these
numbers gives the number of colourings constant on every part of the given
partition.

By Inclusion-Exclusion, we can calculate $|\bigcup_{i\in I}F_i|$ for every
$I\subseteq\{1,\ldots, n\}$ if we know $|\bigcap_{i\in I}F_i|$ for every
such $I$.
\end{pf}

If $m = 2$, $|F_1| = a$ and $|F_2 | = b$, we have a ring of cliques
$R(a,b,1,1)$.
This is a specialisation of a case we have already considered; but, as we saw,
it is general enough to prove the $\alpha + n$ conjecture for all quadratic
integers.

For $m = 3$, we get a six-parameter family of cubics as the ``interesting
factors''. Adam Bohn \cite{r3} has used this family to show that the $\alpha+n$
conjecture is also true for cubic integers.

In general, the ``interesting factor'' has degree $m$ and has $2^m-2$ free
parameters (which must be non-negative integers). Work is proceeding on
using this polynomial to prove further cases of the conjecture. The difficulty
is the exponentially large number of parameters! We hope that this family is
general enough to prove the $\alpha+n$ conjecture.

\section{Other families of graphs}
\label{s5}

In this section we consider some other families of graphs. Unlike the types
considered above (rings of cliques and bicliques), we do not obtain factors of
bounded degree with several free parameters: the parameters appear in the
exponents.

\subsection{Complete bipartite graphs}

The chromatic polynomial of the complete bipartite graph $K_{m,n}$ can be
computed explicitly. Think of $m$ as fixed and $n$ as increasing. Now suppose
that $k$ colours are used on the part of size $m$; the colour classes form a
partition with $k$ parts, and there are $(q-k)^n$ ways to colour the other
part. So the chromatic polynomial is
\[\sum_{k=1}^mS(m,k)(q)_k(q-k)^n = q(q-1)F_{m,n}(q),\]
where $S(m,k)$ is the Stirling number of the second kind (the number of
partitions of $\{1,\ldots,m\}$ into $k$ parts).

For example, we have
\begin{eqnarray*}
F_{2,n}(q) &=& (q-1)^{n-1}+(q-2)^n,\\
F_{3,n}(q) &=& (q-1)^{n-1}+3(q-2)^n+(q-2)(q-3)^n.
\end{eqnarray*}
Note that the degree of the ``interesting'' factor is not bounded by a function
of $m$ in this case.

By computation, we found that, for $3\le n\le 100$, the polynomial $F_{2,n}(q)$
is irreducible if $n$ is not congruent to $2$ mod~$6$; in the remaining cases,
we have $F_{2,n}(q)=F_{2,2}(q)G_n(q)$, where $G_n(q)$ is irreducible. We now
show that at least the factorisation holds in general.

\begin{prop}
$F_{2,2}(q)$ divides $F_{2,6k+2}(q)$ for all $k\ge1$.
\end{prop}

\begin{pf}
Put $x=1-q$. We have $F_{2,2}(q)=(x+1)^2-x=x^2+x+1$, so its roots are primitive
cube roots of unity. If $\omega$ is such a root, then $\omega^3=1$ and
$(\omega+1)^2=\omega$. So, if $n=6k+2$, we have
\[(\omega+1)^n=\omega^{3k+1}=\omega^{6k+1}=\omega^{n-1},\]
so $\omega$ is a root of $F_{2,n}(q)=(x+1)^n-x^{n-1}$. Thus the irreducible
polynomial $F_{2,2}(q)$ divides $F_{2,n}(q)$.
\end{pf}

For $m,n>2$, the polynomial $F_{m,n}(q)$ is irreducible for all the values we
tested.

We note in passing that the Galois group of each of these irreducible
polynomials that we tested is the symmetric group.

\subsection{Theta-graphs}

Let $s$ and $p$ be integers at least $2$. The theta-graph $\Theta^{s,p}$ is
the graph with $2+p(s-1)$ vertices obtained from $p$ disjoint paths of length
$s$ by identifying all the left-hand endpoints and also all the right-hand
endpoints.

These graphs were used in Sokal's proof~\cite{r16} that chromatic roots are
dense in the complex plane. Their chromatic polynomials are known~\cite{r16}; the
chromatic polynomial of $\Theta^{s,p}$ is a product of linear factors and an
``interesting'' factor
\[G_{s,p}(x)=\frac{x(x^s-1)^p-(x^s-x)^p}{x(x-1)^p},\]
where $x=1-q$.

Note that the graph $\Theta^{2,p}$ is just the complete bipartite graph
$K_{2,p}$. Note also that
\[(x-1)G_{s,2}(x) = (x^s-1)^2-x(x^{s-1}-1)^2 = x^{2s-1}-1,\]
so the roots of $G_{s,2}(x)$ are precisely the $(2s-1)$th roots of unity other
than $1$. Indeed, $G_{s,2}(x)$ is the product of the $d$th cyclotomic
polynomials over all $d$ dividing $2s-1$ except for $d=1$.

The result of the preceding section generalises:

\begin{prop}
The polynomial $G_{s,2}(x)$ divides $G_{s,p}(x)$ if and only if $p$ is
congruent to $2$ modulo $2(2s-1)$.
\end{prop}

\begin{pf}
Let $\omega$ be a $(2s-1)$th root of unity other than $1$, and let $p$ be
congruent to $2$ mod $2(2s-1)$. Then $G_{s,2}(\omega) = 0$, so
\[(\omega^s-\omega)^2 = \omega(\omega^s-1)^2,\]
by the calculation preceding the theorem. Let $p = (4s-2)k + 2$, and raise
both sides of this equation to the power $(2s-1)k + 1$. Noting that
$\omega^{2s-1} = 1$, we have
\[(\omega^s − \omega)^p = \omega(\omega^s-1)^p.\]
So every root of $G_{s,2}(x)$ is a root of $G_{s,p}(x)$, and the sufficiency
is proved.

Now let $\omega$ be a primitive $(2s-1)$th root of unity, and suppose that
$G_{s,p}(\omega)=0$. Then
\[(\omega^s-\omega)^p=\omega(\omega^s-1)^p.\]
Suppose that $p=(4s-2)k+2=e$, where $0<e<4s-2$. Then by the first part,
\[(\omega^s-\omega)^{(4s-2)k+2}=\omega(\omega^s-1)^{(4s-2)k+2},\]
so
\[\left(\frac{\omega^s-\omega}{\omega^s-1}\right)^e=1.\]
But $\omega^s-\omega=-\omega^s(\omega^s-1)$, so $\omega^{2se}=1$. This implies
$\omega^e-1$, contrary to assumption.
\end{pf}

The irreducibility of the interesting factor was proved by Delbourgo and
Morgan~\cite{dm-2014}:

\begin{theorem}
Let
\[
P(\Theta^{s,p}; q)=P( \Theta^{a,k+1}, q)=(-1)^{k+1}q (q-1)H_{a}(1-q)
\]
where $s=a$, $p=k+1$ and with a change of variable, $X=1-q$, we have $H_{a}(X) =X^{k+1}-X^{(k+1)a-1}+X^{k}-1$.
Then the interesting factor $G_{a}$, dividing $H_{a}$, is given by the quotient
\[
G_{a} = \cases{%
\frac{(X-1)H_a(X)}{(X^{k+1}-1)(X^d+1)} & if $k$ is odd and $d=\gcd (k-1, 2a-1)>1$,\cr
\frac{H_a(X)}{X^{k+1}-1} & otherwise,\cr
}
\]
and is irreducible over $\mathbb{Q}$.
\end{theorem}

\subsection{Generalised theta graphs}

 We denote the theta graph with consecutive path of lengths 
$ns-n+1, ns-n+2, \ldots, ns$ by $\Theta_{c(s,n)}$ where $a\geq 2$ and
$n\geq 2$.

The chromatic roots of $\Theta_{c(s,n)}$ are closely related to the chromatic 
roots of the theta graph $\Theta^{s,n}$ with  $n$ paths of length $s$.
In \cite{dm-2014}, an explanation for this relationship  is given.
In addition, a description of the Galois group in the case $n=3$ is provided.  
 
After a variable change $x=1-q$, the chromatic polynomial of  $\Theta^{s,n}$ can be expressed as:
\[
P(\Theta^{s,n}, x) = (-1)^{(s+1)n}x(x-1)\times \frac{{f(x)} }{(x-1)^{n}},
\]
where  $f(x)=(x^{s}-1)^{n}-x^{n-1}(x^{s-1}-1)^{n}$.

Similarly, the chromatic polynomial of  $\Theta_{c(s,n)}$ can be expressed
\[
 P(\Theta_{c(s,n)};x)=\frac{P(C_{ns-n+2};x) \ldots P(C_{ns};x)}{P(K_{2};x)^{n-1}}\times (-1)^{ns+1}x\times{g(x)}
\]
where $g(x)=x^{ns}-x^{ns-1}+x^{n-1}-1$ and $C_{i}$ is the cycle of order $i$.

\begin{theorem}
 If $\alpha$ is a root of $g(x)$ then $\alpha^{n}$ is a root of $f(x)$.  \label{thm1}
\end{theorem}

\begin{pf}
\begin{eqnarray*}
 f(\alpha^{n})&=&(\alpha^{ns}-1)^{n}-(\alpha^{n})^{n-1}(\alpha^{n(s-1)}-1)^{n}\\
&=&(\alpha^{ns}-1)^{n}-(\alpha^{ns-1}-\alpha^{n-1})^{n}.
\end{eqnarray*}
As $g(\alpha)=0$ we have $\alpha^{ns}-1=\alpha^{ns-1}-\alpha^{n-1}$ and so
\[
  f(\alpha^{n})=(\alpha^{ns-1}-\alpha^{n-1})^{n}-(\alpha^{ns-1}-\alpha^{n-1})^{n}=0.
\]
\end{pf}

Hence, we have an explanation of the non-trivial relationships between the
chromatic roots of graphs $\Theta^{s,n}$ and chromatic roots of
$\Theta_{c(s,n)}$.   This leads to the following question, a companion to the
$n\alpha$ conjecture: 

\begin{conj}
If $\alpha$ is a chromatic root, then $f(\alpha)=\alpha^{n}$ is a chromatic
root for some $n\in \mathbb{N}$.
\end{conj}

The maxmaxflow $\Lambda$ of a generalised theta graph is the number of
disjoint paths.   
Every chromatic root $q$ of generalised theta graph lies in the disc
$|q-1|\leq \frac{\Lambda-1}{\log 2}=\frac{n-1}{\log 2}$  \cite{rs-2015}.   
It was shown in \cite{bhsw-2001} that $\Theta^{2,n}$  gives the root that
maximises $|q-1|$ over all generalised theta graphs with $n\leq 8$ paths and
conjectured this to be true for larger $n$.   Theorem \ref{thm1} gives some
support to this conjecture, as  it shows that the interesting chromatic roots
of $\Theta^{2,n}$ are larger than the chromatic roots of $\Theta_{c(a,n)}$.  

\section{Galois groups}
\label{s6}

If the $\alpha+n$ conjecture is true, then every transitive permutation group
which actually occurs as a Galois group over the rationals would occur as the
Galois group of an irreducible factor of a chromatic polynomial. The Inverse
Galois Problem asks whether every transitive permutation group actually arises
in this way; we cannot tackle this question, but we have investigated which
small transitive groups arise as Galois groups in the cases we have considered.

For the cases of rings of cliques, and graphs built from families of sets, we
have polynomials of degree bounded in terms of the number of cliques in the
ring, or sets in the family. These cases are amenable to computation. We
have looked at the rings of cliques. Note that the computational technique
we used involved identifying the Galois group as a transitive permutation
group, and is viable for polynomials of degree up to fifteen, that is, for
rings of at most sixteen cliques.

We note that all cyclotomic polynomials will occur here -- the $n$th
cyclotomic polynomial divides the interesting factor of the chromatic
polynomial of an $(n+1)$-cycle. In particular, if $n$ is prime, this
interesting factor is irreducible, with Galois group cyclic of order $n-1$.

The next table shows what happens for small values.

For given $n$, we test all non-decreasing $n$-tuples $(a_1,\ldots,a_n)$ of
positive integers with gcd equal to $1$ and $a_n\le l$. For each such
$n$-tuple, we test the interesting factor of $R(a_1,\ldots,a_n,1)$ (a
polynomial of degree $n-1$). If it is
irreducible, we compute its Galois group. In the table, $S_n$ and $A_n$ are
the symmetric and alternating groups of degree $n$, $C_n$ the cyclic group of
order $n$, $V_4$ the Klein group of order $4$, $D_n$ the dihedral group of
order $2n$, $F_5$ the Frobenius group of order $20$ (the affine group over the
integers mod $5$); and $\wr$ denotes the wreath product of permutation groups.
Entries in the columns labelled ``red'' and ``$S_{n-1}$'' give the number of
tuples for which the polynomial was reducible or had symmetric Galois group;
entries in brackets in the ``Other'' column give these multiplicities for other
groups (if greater than one). The calculations were performed using
GAP \cite{r8}.

The cyclic groups of order $p-1$ (for $p$ prime and $p > 5$) all arise from
the $(p + 2)$-cycle, as explained earlier.

\begin{center}
\begin{tabular}{|l|l|l|l|l|}
\hline
$n$ & $l$ & red & $S_{n-1}$ & Other \\
\hline
$4$ & $30$ & $2581$ & $34471$ & $C_3(\times15)$ \\
$5$ & $30$ & $2677$ & $260658$ & $C_4(\times6)$, $V_4(\times7)$, \\
& & & & $D_4(\times1104)$, $A_4(\times11)$ \\
$6$ & $30$ & $23228$ & $1555851$ & $D_5$, $F_5(\times2)$, $A_5(\times3)$\\
$7$ & $20$ & $2685$ & $642636$ & $C_6$, $S_2\wr S_3(\times10)$,\\
& & & & $S_3\wr S_2(\times145)$, $\mathrm{PGL}(2,5)(\times5)$ \\
$8$ & $10$ & $1132$ & $22630$ & \\
$9$ & $8$ & $152$ & $11054$ & $S_4\wr S_2(\times3)$ \\
$10$ & $8$ & $1061$ & $18089$ & \\
$11$ & $6$ & $29$ & $4248$ & $C_{10}$ \\
$12$ & $6$ & $592$ & $5492$ & \\
$13$ & $6$ & $33$ & $8415$ & $C_{12}$ \\
$14$ & $6$ & $884$ & $10609$ & \\
$15$ & $6$ & $307$ & $15045$ & \\
$16$ & $6$ & $1366$ & $18813$ & \\
\hline
\end{tabular}
\end{center}

Note that we have achieved every transitive permutation group of degree
at most $4$, but for degree $5$ we are missing the cyclic group. The unique
example of the dihedral group $D_5$ occurs for $R(1,4,4,9,9,9,25)$. For degree
$6$, we have seen only five of the $16$ transitive groups.

If certain groups were never realised as Galois groups of chromatic polynomials
of rings of cliques, then this family would not be general enough to prove the
$\alpha+n$ conjecture.

It would be interesting to do similar computation for bicliques.

Further lists of Galois groups of chromatic polynomials can be found in
\cite{r12,dm-2014}.

\subsection{Further speculation}

The Galois group of a ``random'' polynomial is typically the symmetric group
of its degree.

The chromatic polynomial of a random graph cannot be irreducible, since
it will have many linear factors $q-m$, for $m$ up to the chromatic number.
Bollob\'as \cite{r5} showed that the chromatic number is almost surely close
to $n/(2\log_2n)$.

On the basis of admittedly very limited evidence, we propose the following
conjecture:

\begin{conj}
The chromatic polynomial of a random graph is typically a product of linear
factors and one irreducible factor whose Galois group is the symmetric group
of its degree.
\end{conj}

Of course, not all graphs have this property; not even all graphs which
are not clique-separable. There are graphs in which all irreducible factors of
the chromatic polynomial are linear. Chordal graph (graphs in which every
cycle of length greater than 3 has a chord) has this property, and Braun 
\textit{et al.} \cite{r6} conjectured that there are no others; this conjecture
was refuted by Read, who observed that the graph obtained from $K_6$ by
subdividing an edge has chromatic polynomial
\[P_G(q) = q(q-1)(q-2)(q-3)^3(q-4)\]
but is not chordal since it contains an induced $4$-cycle. It seems to be a
difficult open problem to characterize graphs with this property; Dong
\textit{et al.} \cite{r7} give some results for rings of cliques.

On the other hand, Morgan \cite{r12} found that there is a graph on nine
vertices whose chromatic polynomial has two quadratic factors, one with
real roots, and the other with non-real roots. It is labelled $198748$ in the
Geng listing \cite{r10}.

The chromatic polynomial is a specialisation of the two-variable Tutte
polynomial, which itself is a specialisation of the ``multivariate Tutte
polynomial'' which is described in detail in \cite{r17}. This polynomial has a
``local'' variable for each edge of the graph (or more generally, element of
the matroid), and one global variable. It was shown by de Mier \textit{et al.}
\cite{r11} that, for a
connected matroid (in particular, for a $2$-connected graph), the two-variable
Tutte polynomial is irreducible. Furthermore, Bohn \textit{et al.} \cite{r4}
showed that, under the same hypotheses, the multivariate Tutte polynomial
(regarded as a polynomial in the global variable over the field of fractions
of all the local variables) has Galois group the symmetric group. Thus, one
would expect that ``almost all'' specialisations of this polynomial would have
symmetric Galois group. However, we are interested in particular
specialisations, where it is not known whether such a result holds. In
particular, is the Galois group of the two-variable Tutte polynomial of a
connected matroid (as a polynomial in one variable over the field of fractions
of the other) the symmetric group?

\paragraph{Acknowledgment} Much of this work was done in a seminar at the Isaac
Newton Institute in Cambridge, U.K., during a workshop on ``Combinatorics
and Statistical Mechanics'' during the second half of 2008. The participants
in the seminar were Peter Cameron, Vladimir Dokchitser, F. M. Dong, Graham
Farr, Tatiana Gateva-Ivanova, Bill Jackson, Kerri Morgan, Alex Scott,
James Sellers, Alan Sokal, David Wagner, and David Wallace. We are grateful
to the Institute for the excellent facilities and opportunities for interaction
which it provided. The present authors are also grateful to their colleagues
for allowing the results of the seminar to be included in this paper. Our
gratitude to Adam Bohn and Peter M\"uller for helpful comments is also
acknowledged. Finally, we are grateful to the Faculty of Engineering, Computing
and Mathematics at the University of Western Australia for seed funding for a
meeting at which the authors were able to work on a final version of the paper.
This work was supported in part by the Australian Research Council 
under the discovery grant ARC-DP110100957.

\end{document}